\magnification 1200
\baselineskip 14pt
\input amssym.def
\input amssym.tex

\def\qed{\hfill\vbox{\hrule\hbox{\vrule\kern3pt
    \vbox{\kern6pt}\kern3pt\vrule}\hrule}} 

\def\BGS{[BGS]}
\def\Horrocks{[H]}
\def\Knorrer{[K]}
\def\Ottavianiuno{[O1]}
\def\Ottavianidos{[O2]}

\centerline{\bf VECTOR BUNDLES ON $G(1,4)$ WITHOUT
INTERMEDIATE COHOMOLOGY}

\centerline{E. Arrondo and B. Gra\~na}
\bigskip
\bigskip

{\bf \S 0. Introduction}
\bigskip

A known result by Horrocks (see \Horrocks) characterizes the line
bundles on a projective space as the only indecomposable vector
bundles without intermediate cohomology. This result has been
generalized by Ottaviani (see \Ottavianiuno, \Ottavianidos) to
quadrics and Grassmannians. More precisely, he characterizes direct
sums of line bundles as those vector bundles without intermediate
cohomology and satisfying other cohomological conditions.

More generally, Kn\"orrer (see \Knorrer) has proved for any quadric
that the line bundles and spinor bundles (and their twists by line
bundles) are characterized by the property of being indecomposable
and not having intermediate cohomology (there is an unpublished
independent proof of this fact by I. Sols, which has been the
starting point of the present work). Buchweitz, Greuel and Schreyer
(see \BGS) proved a ``converse" of such results: only in the case of
linear spaces and quadrics there are, up to a twist, a finite
number of indecomposable vector bundles without intermediate
cohomology. 

The goal of this paper is to generalize Horrock's result to the
Grassmann variety $G(1,4)$ of lines in ${\Bbb P}^4$, in the sense
of characterizing those vector bundles on it without intermediate
cohomology. 

The paper is distributed in three sections. In the first one, we
give the preliminaries on vector bundles on $G(1,4)$ that will be
needed for the sequel. In the second section, we characterize the
universal bundles on $G(1,4)$ as those without intermediate
cohomology an verifying other cohomology vanishings. Finally, in the
last section we prove our main result, in which --according to the
mentioned result in
\BGS-- we obtain big families of vector bundles without intermediate
cohomology.

We want to acknowledge the tremendous help that has been for us the
Maple package Schubert, created by S.A. Str{\o}mme and S. Katz. Its
use has significantly contributed to an efficient computation of the
cohomology of vector bundles in $G(1,4)$.
\bigskip
\bigskip

{\bf \S 1. Preliminaries}
\bigskip

Let $G=G(1,4)$ denote the Grassmann variety of lines in
${\Bbb P}^4={\Bbb P}(V)$, the projective space of hyperplanes of
$V$. We will assume the ground field to have characteristic
zero, although all our results are likely to hold in any
characteristic different from two. Consider the universal exact
sequence on $G$ defining the universal vector bundles ${\cal Q}$ and
${\cal S}$ of respective ranks two and three:
$$0\to\check{\cal S}\to V\otimes{\cal O}_G\to {\cal Q}\to 0        
             \eqno{(1.1)}$$
\noindent (a check means a dual vector bundle).

The second symmetric power of the above epimorphism induces a long
exact sequence
$$\matrix{
0\to {\cal S}(-1)\to\check{\cal S}\otimes V&&\to&&
S^2V\otimes{\cal O}_G\to S^2Q\to 0\cr
&\searrow&&\nearrow&\cr
&&M&&\cr
&\nearrow&&\searrow&\cr
\hfill0&&&&0\hfill
                   }\eqno{(1.2)}$$
\noindent where the rank-twelve vector bundle $M$ is defined to be the
corresponding kernel, and we made the identification
$\wedge^2\check{\cal S}\cong{\cal S}(-1)$ (as usual we write ${\cal
O}_G(1)\cong\wedge^3{\cal S}\cong\wedge^2{\cal Q}$).

On the other hand, taking the second exterior power in the dual
universal sequence we have the following natural long exact sequence
(defining the rank-seven vector bundle $K$ as a kernel):

$$\matrix{
0\to S^2\check{\cal Q}\to\check{\cal Q} \otimes V^*&&\to&&
\bigwedge^2V^*\otimes{\cal O}_G\to \bigwedge^2{\cal S}\to 0\cr
&\searrow&&\nearrow&\cr
&&K&&\cr
&\nearrow&&\searrow&\cr
\hfill0&&&&0\hfill
                     }\eqno{(1.3)}$$

It is easy to see that ${\rm Ext}^1({\cal S},K)=V^*$, so that,
for any natural numbers $i,j$ there are non-trivial extensions
$$0\to K^{\oplus i}\to {\cal G}\to{\cal S}^{\oplus j}\to
0.\eqno{(1.4)}$$

\noindent{\bf Definition.} An indecomposable direct summand of a
vector bundle ${\cal G}$ as in (1.4) will be called a {\it vector
bundle of type (I)}. 
\bigskip

\noindent {\bf Example 1.1.} Consider the following commutative
diagram of exact sequences coming from (1.3) and the dual of (1.2),
which defines $P$ as a pull-back:

$$\matrix{
 &   &  &   &     0               &   &        0          &   & \cr
 &   &  &   &  \uparrow           &   &  \uparrow         &   & \cr
0&\to&K &\to&\bigwedge^2V^*\otimes{\cal O}_G&\to& \check{\cal S}(1)
&\to&0\cr
 &   &||&   &  \uparrow           &   &    \uparrow       &   & \cr
0&\to&K &\to&     P               &\to&{\cal S}\otimes V^*&\to&0\cr
 &   &  &   &  \uparrow           &   &    \uparrow       &   & \cr
 &   &  &   &  \check M           & = &    \check M       &   & \cr
 &   &  &   &  \uparrow           &   &    \uparrow       &   & \cr
 &   &  &   &     0               &   &        0          &   & 
}$$
\noindent Since Ext$^1({\cal O}_G,\check M)=0$, it follows that
the middle vertical exact sequence splits. This shows that $\check
M$ is a vector bundle of type (I). In fact, the middle horizontal
exact sequence is an element in Ext$^1({\cal S}\otimes V^*,K)\cong
{\rm Hom}(V,V)$, which is represented by the identity map on $V$.

\bigskip
Similarly, one can observe that ${\rm Ext}^1(\check K,K)=V$. This
means that, in general, for a vector bundle ${\cal G}$ appearing in
an exact sequence (split or not) as in (1.4) there are non-trivial
extensions
$$0\to{\cal G}\to{\cal G}'\to\check K^{\oplus l}\to 0.$$

Since ${\rm Ext}^1(\check K,{\cal S})=0$, it follows that such a
${\cal G}'$ appears in an exact sequence

$$0\to K^{\oplus i}\to {\cal G}'\to {\cal S}^{\oplus j}\oplus\check
K^{\oplus l}\to 0.\eqno{(1.5)}$$

\noindent{\bf Definition.} An indecomposable direct summand of a
vector bundle ${\cal G'}$ as in (1.5) will be called a {\it vector
bundle of type (II)}. 
\bigskip

Finally, the left side of the exact sequence (1.3) shows that $K(1)$
is generated by its global sections. This yields the following exact
sequence defining $E(1)$ as a kernel:
$$0\to E(1)\to V\otimes V^*\otimes{\cal O}_G\to K(1)\to
0\eqno{(1.6)}$$
\noindent where $E$ a vector bundle of rank 18. The vector bundle
$E(1)$ has the following non-zero extensions groups:
Ext$^1(K(1),E(1))$ (generated by the extension (1.6)), Ext$^1({\cal
Q},E(1))$, Ext$^1(\check K,E(1))$ and Ext$^1({\cal S},E(1))$. We
give a general definition containing in particular vector bundles
of type (II) and their duals:
\bigskip

\noindent {\bf Definition.} A vector bundle of type (III) will be an
indecomposable direct summand of a (maybe trivial) extension
$$0\to 
E(1)^{\oplus i}\oplus K^{\oplus j}\oplus\check{\cal S}^{\oplus k}
\to {\cal G}\to 
{\cal S}^{\oplus l}\oplus\check K^{\oplus m}\oplus{\cal
Q}^{\oplus n}\to 0$$
\bigskip

\noindent{\bf Example 1.2.} From (1.6) and the first short exact
sequence in (1.3) we obtain the following commutative diagram of
exact sequences defining $P$ as a pull-back:

$$\matrix{
 &   &  &   &     0               &   &        0          &   & \cr
 &   &  &   &  \uparrow           &   &  \uparrow         &   & \cr
0&\to&E &\to&V\otimes V^*\otimes{\cal O}_G(-1)&\to& K     &\to&0\cr
 &   &||&   &  \uparrow           &   &    \uparrow       &   & \cr
0&\to&E &\to&     P      &\to &\check{\cal Q}\otimes V^*   &\to&0\cr
 &   &  &   &  \uparrow           &   &    \uparrow       &   & \cr
 &   &  &   & S^2\check{\cal Q}   & = &S^2\check{\cal Q}  &   & \cr
 &   &  &   &  \uparrow           &   &    \uparrow       &   & \cr
 &   &  &   &     0               &   &        0          &   & 
}$$
Since Ext$^1({\cal O}_G(-1),S^2\check{\cal Q})=0$, it follows that
the middle vertical exact sequence splits, and hence $S^2\check{\cal
Q}$ is of type (III).
\bigskip

\noindent{\bf Definition.} A vector bundle $F$ on $G$ is said
{\it not to have intermediate cohomology} if $H^i(G,F(l))=0$
for all $l\in {\Bbb Z}$ and $i=1,\ldots,5$ (since all cohomology
groups are taken on $G$, we will for short write $H^i(F)$).
\bigskip

\noindent{\bf Remark.} It is easy to see that the vector bundles
${\cal Q}, {\cal S}, \check{\cal S}, K, \check K, E$ are simple
(i.e. their only endomorphisms are multiplications by constants), 
indecomposable and have no intermediate cohomology. This implies 
in particular that all vector bundles of type (III) have no
intermediate cohomology (in fact, all the vector bundles appearing
in this section are, maybe up to a twist, of type (III), as we have
remarked). The goal of this paper is to prove that any vector
bundle on $G$ without intermediate cohomology is obtained in this
way. 
\bigskip

\noindent{\bf Table 1.3.} For the reader's convenience, we list here
the only non-zero intermediate cohomology of the above five vector
bundles when tensored with ${\cal Q}$ and $\check{\cal S}$:

\noindent $h^1({\cal Q}\otimes\check{\cal
S}(-1))=h^5({\cal S}\otimes{\cal Q}(-5))=h^1(\check{\cal
S}\otimes{\cal Q}(-1))=h^2(\check{\cal S}\otimes\check{\cal
S}(-1))=h^2(K\otimes{\cal Q}(-2))=h^3(K\otimes\check{\cal
S}(-2)=h^4(\check K\otimes{\cal Q}(-4))=h^5(\check
K\otimes\check{\cal S}(-4))=h^3(E\otimes{\cal
Q}(-2))=h^4(E\otimes\check{\cal S}(-2)=1.$

\noindent $h^1(K\otimes\check{\cal S})=h^1(E\otimes{\cal
Q})=h^2(E\otimes\check{\cal S})=h^1(E\otimes\check{\cal S}(1))=5.$

Most of the above equalities can be derive from the others by using
the universal exact sequence (1.1) or the Serre duality, taking into
account that the canonical line bundle on $G$ is $\omega_G={\cal
O}_G(-5)$.

\bigskip
\bigskip

\centerline{\bf \S2. Characterization of the universal bundles}
\bigskip

We start by recalling Ottaviani's characterization of direct sums
of line bundles, when particularized to $G(1,4)$.

\proclaim Theorem 2.1. (Ottaviani, \Ottavianiuno, \Ottavianidos)
Let $F$ be a vector bundle on $G(1,4)$. Then $F$ is a direct sum of
line bundles if and only if the following conditions hold:
\item{a)} $F$ has no intermediate cohomology
\item{b)} $H^i(F\otimes{\cal Q}(l))=0$ for any $i=1,\ldots,5$ and
$l\in{\Bbb Z}$
\item{c)} $H^1(F\otimes\check{\cal S}(l))=0$ for any $l\in{\Bbb Z}$.

\bigskip

\noindent {\bf Remark.} Ottaviani's original statement is not as we
gave it. Instead of conditions b) and c), his conditions are (see
\Ottavianiuno, Theor. 1 (c) for $k=1$, $n=4$):
\itemitem{b')} $H^i(F\otimes{\cal S}(l))=0$ for any $i=3,4,5$
and $l\in{\Bbb Z}$
\itemitem{c')} $H^i(F\otimes\check{\cal S}(l))=0$ for any
$i=1,2,3$ and $l\in{\Bbb Z}$.

\noindent These are clearly equivalent to b) and c) in our statement
by taking cohomology in the universal exact sequence (1.1) and its
dual tensored with $F(l)$, and using the assumption that $F$ has no
intermediate cohomology. 

The idea for proving the main theorem is to successively remove the
six extra cohomological conditions appearing in b) and c) to
eventually characterize those vector bundles without intermediate
cohomology. Each time we remove a condition, a new family of vector
bundles will appear. Table 1.3 indicates which vector bundles must
appear each time. We will characterize in this section the universal
vector bundles, by removing --one by one-- the conditions in Theorem
2.1 that they do not satisfy. 

The first condition we remove will be c), which is the only
condition that ${\cal Q}$ does not verify. Hence, ${\cal Q}$ should
be  characterized by conditions a) and b) in Ottaviani's theorem.
Notice that then we obtain a result completely analogous to Horrock's
theorem, the role of line bundles being played now by line
bundles and their tensor products with
${\cal Q}$. We will prove this result in detail, the others being
sketched as long as they are similar (in fact, this proof will have
a difficulty at the beginning not appearing in the remaining proofs).
The precise statement is:

\proclaim Theorem 2.2. Let $F$ be an indecomposable vector bundle on
$G(1,4)$. Then
$F$ is, up to a twist with a line bundle, either the trivial line
bundle or ${\cal Q}$ if and only if the
following conditions hold:
\item{a)} $F$ has no intermediate cohomology
\item{b)} $H^i(F\otimes{\cal Q}(l))=0$ for any $i=1,\ldots,5$ and
$l\in{\Bbb Z}$.

\noindent{\it Proof.} Let $F$ be a vector bundle satisfying a) and
b). We will prove our result by induction on $\sum_l
h^1(F\otimes\check{\cal S}(l))$. If this sum is zero, then we are in
the hypotheses of Ottaviani's Theorem 2.1, so that $F$ is a
direct sum of line bundles.

So assume that $h^1(F\otimes\check{\cal S}(l))\neq 0$ for some
$l\in{\Bbb Z}$. By changing if necessary $F$ with a twist, we can
assume $l=0$. Then we have a non-zero element in ${\rm Ext}^1({\cal
S},F)$, which yields a non-trivial extension
$$0\to F\to P\to{\cal S}\to 0.\eqno{(2.1)}$$
We first claim that $P$ also verifies conditions a) and b).
Indeed the only vanishing to check is that of $H^5(P\otimes{\cal
Q}(-5))$, since $h^5({\cal S}\otimes{\cal Q}(-5))=1$. To prove this
vanishing, we first dualize and suitably twist the exact sequences
(1.1), (1.2), (1.3), and using the natural isomorphisms $\bigwedge^2
{\cal S}\cong \check{\cal S}(1)$ and $S^2\check{\cal Q}\cong S^2{\cal
Q}(-2)$ we get a long exact sequence
$$
0\to{\cal Q}(-5)\to{\cal O}_G(-4)^{\oplus 5}\to
 {\cal O}_G(-3)^{\oplus 10}\to{\cal Q}(-3)^{\oplus 5}\to
  {\cal O}_G(-1)^{\oplus 15}\to
  \check M(-1)\to 0.
$$
\noindent The fact that $F$ and ${\cal S}$ (and hence also $P$)
satisfy conditions a) and, for $i=2,3$, also b) easily implies that
there is a commutative diagram 
$$\matrix{
H^5(P\otimes{\cal Q}(-5))&\to& H^5({\cal S}\otimes{\cal Q}(-5))\cr
\uparrow&&\uparrow\cr
H^1(P\otimes\check M(-1))&\to& H^1({\cal S}\otimes\check M(-1))
}$$
\noindent where the vertical maps are isomorphisms. Since we just
need to prove that the map on the top is zero, it suffices to prove
the same for the map in the bottom. By looking at the dual of the
first short exact sequence in (1.2), that map appears in a
commutative diagram
$$\matrix{
H^0(P\otimes\check{\cal S})&\to& H^0({\cal S}\otimes\check{\cal
S})\cr
\downarrow&&\downarrow\cr
H^1(P\otimes\check M(-1))&\to& H^1({\cal S}\otimes\check M(-1)).
}$$
\noindent The claim follows by observing that the vertical map on the
left is an epimorphism (its cokernel lies in $H^1(P\otimes{\cal
S}(-1)\otimes V^*)=H^2(P\otimes\check{\cal Q}(-1)\otimes V^*)=0$),
while the map on the top is zero since the extension (2.1) was
non-trivial and ${\cal S}$ is simple.

It is also immediate to check that $h^1(P\otimes\check{\cal S}(l))=
h^1(F\otimes\check{\cal S}(l))$ for any $l$ except $l=0$, for which
$h^1(P\otimes\check{\cal S})=h^1(F\otimes\check{\cal S})-1$. The
latter follows from the exact sequence
$$
H^0(P\otimes\check{\cal S})\to
H^0({\cal S}\otimes\check{\cal S})\to H^1(F\otimes\check{\cal S})\to
H^1(P\otimes\check{\cal S})\to H^1({\cal S}\otimes\check{\cal S})=0,
$$
\noindent in which, as we observed, the first map is zero and
$h^0({\cal S}\otimes\check{\cal S})=1$.

We can therefore apply the induction hypothesis to $P$ and conclude
that it decomposes as a direct sum of summands of the type ${\cal
O}_G(l)$ and ${\cal Q}(l)$. We next consider the following
commutative diagram defining $P'$ as a pull-back (and in which the
right vertical map is the dual of (1.1)):
$$\matrix{
 &   &  &   &     0        &   &        0            &   & \cr
 &   &  &   &  \uparrow    &   &    \uparrow         &   & \cr
0&\to&F &\to&     P        &\to&    {\cal S}         &\to&0\cr
 &   &||&   &  \uparrow    &   &    \uparrow         &   & \cr
0&\to&F &\to&     P'       &\to&{\cal O}_G^{\oplus 5}&\to&0\cr
 &   &  &   &  \uparrow    &   &    \uparrow         &   & \cr
 &   &  &   &\check{\cal Q}& = &   \check{\cal Q}    &   & \cr
 &   &  &   &  \uparrow    &   &    \uparrow         &   & \cr
 &   &  &   &     0        &   &        0            &   & 
}$$

The middle horizontal exact sequence splits since $F$ has not
intermediate cohomology. The middle vertical exact sequence
also splits, since Ext$^1(P,\check{\cal Q})\cong H^5(P\otimes
{\cal Q}(-5))^*$=0. Hence $P\oplus\check{\cal Q}\cong F\oplus{\cal
O}_G^{\oplus 5}$, from which the theorem follows.
\qed

\bigskip

\proclaim Theorem 2.3. Let $F$ be an indecomposable vector bundle on
$G(1,4)$. Then
$F$ is, up to a twist with a line bundle, either the trivial line
bundle, or ${\cal Q}$, or ${\cal S}$ if and only if the
following conditions hold:
\item{a)} $F$ has no intermediate cohomology
\item{b)} $H^i(F\otimes{\cal Q}(l))=0$ for any $i=1,2,3,4$ and
$l\in{\Bbb Z}$.
\bigskip

{\it Proof.} We prove it by induction on $\sum_l h^5(F\otimes{\cal
Q}(l))$, the zero case being Theorem 2.2. Hence we assume 
$h^5(F\otimes{\cal Q}(l))\neq 0$ for some $l$, and we can suppose
without loss of generality that $l=-5$. Hence, by Serre duality,
there is a non-zero element in Ext$^1({\cal Q},\check F)$. This
produces a non-trivial extension
$$0\to \check F\to P\to{\cal Q}\to 0.\eqno{(2.2)}$$
\noindent It is now immediate (in contrast with the proof of Theorem
2.2) to prove that $\check P$ still satisfies conditions a) and b)
(this is because ${\cal Q}\otimes{\cal Q}$ has no intermediate
cohomology). Also, it holds that
$h^5(\check P\otimes{\cal Q}(l))=h^5(F\otimes{\cal Q}(l))$ for any
$l$, except for $l=-5$ for which
$h^5(\check P\otimes{\cal Q}(-5))=h^5(F\otimes{\cal Q}(-5))-1$ (the
proof being, by using Serre duality, as in Theorem 2.2). Therefore,
by induction hypothesis, $\check P$ decomposes as a direct sum of
summands of the type ${\cal O}_G(l)$, ${\cal Q}(l)$ and ${\cal
S}(l)$.

We consider now the following commutative diagram defining $P'$ as a
pull-back:

$$\matrix{
 &   &        &   &     0        &   &        0            &   & \cr
 &   &        &   &  \uparrow    &   &    \uparrow         &   & \cr
0&\to&\check F&\to&     P        &\to&    {\cal Q}         &\to&0\cr
 &   &   ||   &   &  \uparrow    &   &    \uparrow         &   & \cr
0&\to&\check F&\to&     P'       &\to&{\cal O}_G^{\oplus 5}&\to&0\cr
 &   &        &   &  \uparrow    &   &    \uparrow         &   & \cr
 &   &        &   &\check{\cal S}& = &   \check{\cal S}    &   & \cr
 &   &        &   &  \uparrow    &   &    \uparrow         &   & \cr
 &   &        &   &     0        &   &        0            &   & 
}$$

As in the proof of Theorem 2.2, the middle horizontal exact sequence
splits. Therefore, our result will follow if the middle vertical
exact sequence also splits (this is the only difficulty that did
not appear in Theorem 2.2 and that will appear in the rest of the
proofs). To see this, we study the direct summands of
Ext$^1(P,\check{\cal S})\cong H^1(\check P\otimes
\check{\cal S})$ corresponding to the decomposition of $P$ into
direct summands. Only a summand of the type ${\cal Q}\subset P$ (if
it exists) produces a non-zero summand Ext$^1({\cal Q},\check
{\cal S})\subset{\rm Ext}^1(P,\check{\cal S})$. But then the
corresponding component of the element of $\xi\in{\rm
Ext}^1(P,\check{\cal S})$ defined by the vertical extension must be
zero. Indeed, $\xi$ is the image of the universal extension  (1.1)
under the map
$\pi^*:{\rm Ext}^1({\cal Q},\check{\cal S})\to {\rm
Ext}^1(P,\check{\cal S})$ induced by the projection
$\pi:P\to {\cal Q}$ in (2.2). Since (2.2) is non-split and the only
endomorphisms of ${\cal Q}$ are the multiplications by a constant, it
follows that the restriction of $\pi$ to any direct summand ${\cal
Q}\subset P$ is zero.

\qed

\bigskip
Finally, we prove a theorem characterizing the universal vector
bundles on $G(1,4)$ by means of their cohomology vanishings.

\proclaim Theorem 2.4. Let $F$ be an indecomposable vector bundle on
$G(1,4)$. Then $F$ is, up to a twist with a line bundle, either the
trivial line bundle, or ${\cal Q}$, or ${\cal S}$, or $\check {\cal
S}$ if and only if the following conditions hold:
\item{a)} $F$ has no intermediate cohomology
\item{b)} $H^i(F\otimes{\cal Q}(l))=0$ for any $i=2,3,4$ and
$l\in{\Bbb Z}$.
\bigskip

{\it Proof.} We use now induction on $\sum_l h^1(F\otimes{\cal
Q}(l))$, the zero case being now Theorem 2.3. Again we can assume 
$h^1(F\otimes{\cal Q}(-1))\neq 0$. Therefore, there is a non-zero
element in Ext$^1({\cal Q},F)$, which yields a non-trivial extension
$$0\to F\to P\to{\cal Q}\to 0.$$
\noindent As before, $P$ still satisfies conditions a) and b) and
$h^1(P\otimes{\cal Q}(l))=h^1(F\otimes{\cal Q}(l))$
for any $l$, except for $l=-1$ for which
$h^1(P\otimes{\cal Q}(-1))=h^1(F\otimes{\cal Q}(-1))-1$. Hence,
by induction hypothesis, $P$ decomposes as a direct sum of summands
of the type ${\cal O}_G(l)$, ${\cal Q}(l)$, ${\cal S}(l)$ and
$\check{\cal S}(l)$.

Finally, we consider the following commutative diagram defining
$P'$ as a pull-back:

$$\matrix{
 &   &        &   &     0        &   &        0            &   & \cr
 &   &        &   &  \uparrow    &   &    \uparrow         &   & \cr
0&\to&    F   &\to&     P        &\to&    {\cal Q}         &\to&0\cr
 &   &   ||   &   &  \uparrow    &   &    \uparrow         &   & \cr
0&\to&    F   &\to&     P'       &\to&{\cal O}_G^{\oplus 5}&\to&0\cr
 &   &        &   &  \uparrow    &   &    \uparrow         &   & \cr
 &   &        &   &\check{\cal S}& = &   \check{\cal S}    &   & \cr
 &   &        &   &  \uparrow    &   &    \uparrow         &   & \cr
 &   &        &   &     0        &   &        0            &   & 
}$$

The middle horizontal exact sequence splits as usual, and the
splitting of the middle vertical exact sequence is proved as in
Theorem 2.3. This proves the theorem.
\qed

\bigskip
\bigskip

\centerline{\bf \S3. Vector bundles without intermediate cohomology}
\bigskip

We continue here the strategy of the previous section. The
difference now is that we will not obtain a finite number of vector
bundles when removing any of the conditions in Theorem 2.4 b).

\bigskip
\proclaim Theorem 3.1. Let $F$ be an indecomposable vector bundle on
$G(1,4)$. Then
$F$ is, up to a twist with a line bundle, either the trivial line
bundle or ${\cal Q}$, ${\cal S}$, $\check {\cal S}$ or a vector
bundle of type (I) if and only if the following conditions hold:
\item{a)} $F$ has no intermediate cohomology
\item{b)} $H^i(F\otimes{\cal Q}(l))=0$ for any $i=3,4$ and
$l\in{\Bbb Z}$.
\bigskip

{\it Proof.} We use induction on $\sum_l h^2(F\otimes{\cal
Q}(l))$, the zero case being Theorem 2.4. We can assume 
$h^2(F\otimes{\cal Q}(-2))\neq 0$. Hence, by using the universal
exact sequence (1.1) there is a non-zero element in
Ext$^1(\check{\cal S}(1),F)$, which yields a non-trivial extension
$$0\to F\to P\to\check{\cal S}(1)\to 0.$$
\noindent As in the proofs of the preceeding section, $P$ still
satisfies conditions a) and b) and
$h^2(P\otimes{\cal Q}(l))=h^2(F\otimes{\cal Q}(l))$
for any $l$, except for $l=-2$ for which
$h^2(P\otimes{\cal Q}(-2))=h^2(F\otimes{\cal Q}(-2))-1$. By induction
hypothesis, $P$ decomposes as a direct sum of summands of the type
${\cal O}_G(l)$, ${\cal Q}(l)$, ${\cal S}(l)$, $\check{\cal S}(l)$
and twists of vector bundles of type (I).

Now we consider the following commutative diagram defining
$P'$ as a pull-back, and in which the right column is the second
short exact sequence in (1.3):

$$\matrix{
 &   &      &   &     0        &   &        0             &   & \cr
 &   &      &   &  \uparrow    &   &    \uparrow          &   & \cr
0&\to&  F   &\to&     P        &\to&\check{\cal S}(1)     &\to&0\cr
 &   & ||   &   &  \uparrow    &   &    \uparrow          &   & \cr
0&\to&  F   &\to&     P'       &\to&{\cal O}_G^{\oplus 10}&\to&0\cr
 &   &      &   &  \uparrow    &   &    \uparrow          &   & \cr
 &   &      &   &     K        & = &        K             &   & \cr
 &   &      &   &  \uparrow    &   &    \uparrow          &   & \cr
 &   &      &   &     0        &   &        0             &   & 
}$$

The middle horizontal exact sequence splits as usual, because $F$ has
not intermediate cohomology. As for the middle vertical exact
sequence, the main difference now with the proofs in the preceeding
sections is that the element Ext$^1(P,K)$ defining the extension can
have non-zero components diferent from those corresponding to
possible summands $\check{\cal S}(1)\subset P$ (which as usual we
know to produce a zero coordinate). Indeed, any direct summand of
type (I) (we include here possible direct summands ${\cal S}\subset
P$) will yield a non-zero summand of Ext$^1(P,K)$. We decompose
$P=P_1\oplus P_2$, where $P_1$ is a sum of vector bundles of type
(I) and $P_2$ does not have any summand of type (I). Since our
extension lives in Ext$^1(P,K)={\rm Ext}^1(P_1,K)\oplus {\rm
Ext}^1(P_2,K)$ and the second coordinate is zero, it follows that
$P'\cong P'_1\oplus P_2$, where $P'_1$ appears in an exact sequence
$$0\to K\to P'_1\to P_1\to 0.$$  
\noindent But from this exact sequence it is immediate to see that
$P'_1$ is a direct sum of vector bundles of type (I) (since $P_1$ is
too). This completes the proof.
\qed

\bigskip

\proclaim Theorem 3.2. Let $F$ be an indecomposable vector bundle on
$G(1,4)$. Then
$F$ is, up to a twist with a line bundle, either the trivial line
bundle, or ${\cal Q}$, or ${\cal S}$, or $\check {\cal S}$, or a
vector bundle of type (II), or the dual of a vector bundle of type
(II) if and only if the following conditions hold:
\item{a)} $F$ has no intermediate cohomology
\item{b)} $H^3(F\otimes{\cal Q}(l))=0$ for any $l\in{\Bbb Z}$.
\bigskip

{\it Proof.} We use induction on
$\sum_l h^4(F\otimes{\cal Q}(l))$, the zero case being Theorem 3.1.
We can assume 
$h^2(\check F\otimes\check{\cal Q}(-1))=h^4(F\otimes{\cal Q}(-4))\neq
0$. Hence, by using the universal exact sequence (1.1) there is a
non-zero element in Ext$^1(\check{\cal S}(1),\check F)$, which yields
a non-trivial extension
$$0\to \check F\to P\to\check{\cal S}(1)\to 0.$$
\noindent As usual, $\check P$ still
satisfies conditions a) and b) and
$h^4(\check P\otimes{\cal Q}(l))=h^4(F\otimes{\cal Q}(l))$
for any $l$, except for $l=-4$ for which
$h^4(P\otimes{\cal Q}(-4))=h^4(F\otimes{\cal Q}(-4))-1$. By induction
hypothesis, $P$ decomposes as a direct sum of summands of the type
${\cal O}_G(l)$, ${\cal Q}(l)$, ${\cal S}(l)$, $\check{\cal S}(l)$
and twists of vector bundles of type (II) or their duals.

Now we consider the following commutative diagram defining
$P'$ as a pull-back:

$$\matrix{
 &   &        &   &     0        &   &        0             &   & \cr
 &   &        &   &  \uparrow    &   &    \uparrow          &   & \cr
0&\to&\check F&\to&     P        &\to&\check{\cal S}(1)     &\to&0\cr
 &   &   ||   &   &  \uparrow    &   &    \uparrow          &   & \cr
0&\to&\check F&\to&     P'       &\to&{\cal O}_G^{\oplus 10}&\to&0\cr
 &   &        &   &  \uparrow    &   &    \uparrow          &   & \cr
 &   &        &   &     K        & = &        K             &   & \cr
 &   &        &   &  \uparrow    &   &    \uparrow          &   & \cr
 &   &        &   &     0        &   &        0             &   & 
}$$

Once more, the middle horizontal exact sequence splits, and the only
non-zero components (corresponding to direct summands of $P$) in the
element Ext$^1(P,K)$ defining the middle vertical extension can be
those coming from possible summands of type (II) of $P$ (as in the
preceeding proofs, the coordinates corresponding to possible factors
$\check{\cal S}(1)$ are zero). As in the proof of Theorem 3.1, we
also consider ${\cal S}$ and $\check K$ as vector bundles of type
(II). We decompose $P=P_1\oplus P_2$, where $P_1$ is a sum of vector
bundles of type (II) and $P_2$ does not have any summand of type
(II). Hence, it follows that $P'\cong P'_1\oplus P_2$, where $P'_1$
appears in an exact sequence
$$0\to K\to P'_1\to P_1\to 0.$$  
\noindent But, again as in the proof of Theorem 3.1, $P'_1$ is a
direct sum of vector bundles of type (II), which completes the proof.
\qed
\bigskip

We can finally state and prove our result characterizing vector
bundles on $G(1,4)$ without intermediate cohomology.
\bigskip

\proclaim Theorem 3.3. An indecomposable vector bundle on
$G(1,4)$ without intermediate cohomology is, up to a twist with a
line bundle, a vector bundle of type (III).
\bigskip

\noindent{\it Proof.} Let $F$ be a vector bundle on $G(1,4)$ without
intermediate cohomology. We use induction on
$\sum_l h^3(F\otimes{\cal Q}(l))$, the zero case being Theorem 3.2.
Without loss of generality, we assume $h^3(F\otimes{\cal Q}(-3))\neq
0$. Hence, by using the dual of the exact sequences (1.1) and
(1.3) and the identification $\wedge^2{\cal S}\cong \check S(1)$,
there is a non-zero element in Ext$^1(K(1),F)$, which yields
a non-trivial extension
$$0\to F\to P\to K(1)\to 0.\eqno{(3.1)}$$
\noindent As usual, $\check P$ has no intermediate cohomology and
$h^3(\check P\otimes{\cal Q}(l))=h^3(F\otimes{\cal Q}(l))$
for any $l$, except for $l=-3$ for which
$h^3(P\otimes{\cal Q}(-3))=h^3(F\otimes{\cal Q}(-3))-1$. Hence, by
induction hypothesis, $P$ decomposes as a direct sum of summands of
twists of vector bundles of type (III).

Now we consider the following commutative diagram defining $P'$ as a
pull-back, and in which the right column is the exact sequence (1.6):

$$\matrix{
 &   &        &   &     0        &   &        0             &   & \cr
 &   &        &   &  \uparrow    &   &    \uparrow          &   & \cr
0&\to&   F    &\to&     P        &\to&         K(1)         &\to&0\cr
 &   &   ||   &   &  \uparrow    &   &    \uparrow          &   & \cr
0&\to&   F    &\to&     P'       &\to&{\cal O}_G^{\oplus 25}&\to&0\cr
 &   &        &   &  \uparrow    &   &    \uparrow          &   & \cr
 &   &        &   &    E(1)      & = &        E(1)          &   & \cr
 &   &        &   &  \uparrow    &   &    \uparrow          &   & \cr
 &   &        &   &     0        &   &        0             &   & 
}$$

The middle horizontal exact sequences splits once more because $F$
has not intermediate cohomology. As for the middle vertical exact
sequence, we decompose $P=P_1\oplus P_2$, where
$P_1$ is a sum of vector bundles of type (III) and the summands
of $P_2$ are not of type (III), but twists of vector
bundles of type (III) with a non-trivial line bundle. As in the other
proofs in this section, it suffices to prove that the element
$\xi\in{\rm Ext}^1(P,E(1))={\rm Ext}^1(P_1,E(1))\oplus{\rm
Ext}^1(P_2,E(1))$ corresponding to the middle vertical extension
has zero as its second component in this decomposition. 

The extra difficulty now is that, if a summand $H\subset P_2$
produces a non-zero summand ${\rm Ext}^1(H,E(1))\subset{\rm
Ext}^1(P_2,E(1))$, it is not neccesarily $H=K(1)$, but it could
happen, more generally, that $H(-1)$ is a vector bundle of type
(III). In this case, $H$ is a direct summand of a vector bundle
${\cal G}(1)$ fitting in an exact sequence
$$0\to E(2)^{\oplus i}\oplus K(1)^{\oplus j}\oplus\check{\cal
S}(1)^{\oplus k}\to{\cal G}(1)\to{\cal S}(1)^{\oplus l}\oplus
\check K(1)^{\oplus m}\oplus{\cal Q}(1)^{\oplus n}\to0$$

By contradiction, assume that the component of $\xi$ in
$Ext^1(H,E(1))$ is not zero. Then, at least one summand $K(1)$ of
that exact sequence must produce a non-zero map
$$K(1)\subset E(2)^{\oplus i}\oplus K(1)^{\oplus j}\oplus\check{\cal
S}(1)^{\oplus k}\to {\cal G}(1)\to H\subset P\to K(1)$$
\noindent (the projection $P\to K(1)$ being that of (3.1)). Since the
only endomorphisms of $K$ are multiplications by constants, this
implies that the projection $P\to K(1)$ has a section, so that
the extension (3.1) must be trivial, which is a contradiction.
\qed

\bigskip
\bigskip
{\bf References}
\bigskip

\item{\BGS} Buchweitz, R.O. -- Greuel, G.M. -- Schreyer, F.O., {\it
Cohen-Macaulay modules on hypersurface singularities II}, Invent.
Math. {\bf 88} (1987), 165-182.
\item{\Horrocks} Horrocks, G., {\it Vector bundles on the punctured
spectrum of a ring}, Proc. London Math. Soc. (3) {\bf 14} (1964),
689-713.
\item{\Knorrer} Kn\"orrer, H., {\it Cohen-Macaulay modules on
hypersurface singularities I}, Invent. Math. {\bf 88} (1987),
153-164.
\item{\Ottavianiuno} Ottaviani, G., {\it Crit\`eres de scindage pour
les fibr\'es vectoriels sur les grassmannianes et les quadriques},
C.R. Acad. Sci. Paris, t. {\bf 305}, S\'erie I (1987), 257-260.
\item{\Ottavianidos} Ottaviani, G. {\it Some extensions of Horrocks
criterion to vector bundles on Grassmannians and quadrics}, Annali
Mat. Pura Appl. (IV) {\bf 155} (1989), 317-341.
\bigskip
\bigskip

\centerline{Departamento de Algebra}
\centerline{Facultad de Ciencias Matem\'aticas}
\centerline{Universidad Complutense de Madrid} 
\centerline{28040 Madrid, Spain}
\centerline{\tt enrique@sunal1.mat.ucm.es \ \ \ \
beagra@sunal1.mat.ucm.es}

\end